\documentclass[preprint,11pt]{elsarticle}

\usepackage{amssymb}
\usepackage{amsmath}
\usepackage{amsthm}

\begin{document}

\def\m{{\mathfrak m}} 

\def\p{{\mathfrak p}}    
\def\q{{\mathfrak q}}    
\def\n{{\mathfrak n}} 

\def\S{{\mathcal S}} 
\def\T{{\mathcal T}} 
\def\U{{\mathcal U}} 
\def\O{{\mathcal O}} 
\def\R{{\mathcal R}} 
\def\E{{\mathcal E}} 
\def\B{{\mathcal B}}

\def\I{{\mathcal I}} 
\def\v{{\mathfrak v}}
\def\V{{\mathcal V}}
\def\W{{\mathcal w}}

\def\N{\mathbb N} 

\def\SN{\mathbb N} 
\def\SQ{{\mathbb Q}} 
\def\SZ{\mathbb Z} 
\def\SF{\mathbb F} 
\def\SR{\mathbb R} 
\def\SC{\mathbb C} 

\def\spec{\operatorname{Spec }} 

\def\ff{\frak}
\def\Spec{\mbox{\rm Spec }}
\def\Proj{\mbox{\rm Proj }}
\def\Rees{\mbox{\rm Rees }}

\def\hgt{\mbox{\rm ht }}
\def\type{\mbox{ type}}
\def\Hom{\mbox{ Hom}}
\def\rank{\mbox{ rank}}
\def\Ext{\mbox{ Ext}}
\def\Ker{\mbox{ Ker}}
\def\Max{\mbox{\rm Max}}
\def\End{\mbox{\rm End}}
\def\xpd{\mbox{\rm xpd}}
\def\Ass{\mbox{\rm Ass}}
\def\emdim{\mbox{\rm emdim}}
\def\epd{\mbox{\rm epd}}
\def\repd{\mbox{\rm rpd}}
\def\ord{\mbox{\rm ord}}

\newtheorem{theorem}{Theorem}[section]
\newtheorem{lemma}[theorem]{Lemma}
\newtheorem{proposition}[theorem]{Proposition}
\newtheorem{corollary}[theorem]{Corollary}
\newtheorem{problem}[theorem]{Problem}

\newtheorem{defi}[theorem]{Definitions}
\newtheorem{definition}[theorem]{Definition}
\newtheorem{remark}[theorem]{Remark}
\newtheorem{example}[theorem]{Example}
\newtheorem{question}[theorem]{Question}
\newtheorem{comment}[theorem]{Comment}
\newtheorem{comments}[theorem]{Comments}
\newtheorem{construction}[theorem]{Construction}
\newtheorem{notation}[theorem]{Notation}
\newtheorem{facts}[theorem]{Facts}

\newtheorem{discussion}[theorem]{Discussion}

\newproof{pf}{Proof}

\renewcommand{\thedefi}{}

\begin{frontmatter}



\long\def\alert#1{\smallskip{\hskip\parindent\vrule%
\vbox{\advance\hsize-2\parindent\hrule\smallskip\parindent.4\parindent%
\narrower\noindent#1\smallskip\hrule}\vrule\hfill}\smallskip}

\title{Noetherian intersections of regular local rings of dimension two}

\author[Heinzer]{William Heinzer}
\address[Heinzer]{Department of Mathematics, Purdue University, West
Lafayette, Indiana 47907-1395 U.S.A.}
\ead{heinzer@purdue.edu}

\author[Olberding]{Bruce Olberding} 
\address[Olberding]{Department of Mathematical Sciences, New Mexico State University,
 Las Cruces, NM 88003-8001 U.S.A.}
 \ead{olberdin@nmsu.edu}




\begin{abstract}
  Let $D$ be a 2-dimensional regular local ring and let $Q(D)$ denote the quadratic tree of 2-dimensional 
  regular local overrings of $D$.  We examine the  Noetherian rings that are intersections of rings in $Q(D)$. To do so, we describe the desingularization of projective models over $D$ both algebraically in terms of the saturation of complete ideals and order-theoretically in terms of the quadratic tree $Q(D)$. 

\end{abstract}


\end{frontmatter}


  \section{Introduction}

  Let $D$ be a 2-dimensional regular local ring with quotient field $F$.    
  This article concerns the structure of the Noetherian rings   that are intersections of $2$-dimensional 
  regular local  rings between $D$ and $F$.
    As an intersection of normal rings, such rings  are necessarily normal. 
    We  show these rings  have the property that every maximal ideal has  height~2.  
   Conversely,    it follows from Lipman's work \cite{L1969} on rational singularities that 
     every  normal Noetherian overring $R$ of $D$ with height~2 maximal ideals
     has the form $R = \bigcap_{T \in \U}T$,  where $\U$ is a set of 
      2-dimensional regular local overrings  of $D$. 
      
 In a  paper in preparation \cite{HLO},  we show the existence of subsets $\U$ of the set of  2-dimensional 
 regular local overrings of $D$ such that $R = \bigcap_{T \in \U}T$ is not Noetherian. 
 The question arises as to which sets $\U$ correspond  to normal Noetherian overrings of $D$.
  This question is the main focus of the article.  To address it we situate the problem 
    in the context of the quadratic tree of $D$, that is,  the partially ordered set $Q(D)$ of $2$-dimensional regular local rings that birationally dominate  $D$. 
    
    Each ring $R$ in $Q(D)$  is obtained via Abyhankar-Zariski factorization by a sequence $D = D_0 \subseteq D_1 \subseteq \cdots \subseteq D_n = R$ of quadratic transformations of $D$. The set $\{D_0,\ldots, D_n\}$ is precisely the set of $2$-dimensional regular local rings between $D$ and $R$. In Section 3 we recall basic properties of the tree $Q(D)$, and in Section 4 we recall the concept  of a projective model over $D$.  A nonsingular projective model $X = \Proj D[It]$, with $I$ a complete ${\ff m}_D$-primary ideal,         
        can  be expressed (via its closed points) in terms of $Q(D)$.   We do this in   
        Theorem~\ref{cor5.2} using  the Rees valuation rings and base points of 
        the ideal   $I$.  
    
        Nonsingular projective models over $D$  are central to our approach for describing the
  sets $\U$ in $Q(D)$ that give rise to Noetherian rings $R = \bigcap_{T \in \U}T$. We recall in Proposition~\ref{sat = ns} that every normal projective model $X = \Proj D[It]$ over $D$ has a desingularization, and in Theorem~\ref{thm2.6} we
  use Zariski's structure theorem for complete ideals to
   describe how to obtain the unique minimal desingularization  via saturation of the ideal  $I$.   This leads in Theorem~\ref{order-theoretic} to a strictly order-theoretic description of the closed points of the minimal desingularization of $X$ in terms of the partially ordered set $Q(D)$.   
  

Theorem~\ref{prop3.7}
describes  properties of  the intersection of the $2$-dimensional regular 
         local rings in  an affine component of $X$.  
         Corollary~\ref{thm3.1}  asserts the following description   
         of an intersection of finitely many rings in $Q(D)$:   If   $n$ is  a positive integer and 
  $R$ is an irredundant intersection of $n$ 
  elements in $Q(D)$, then
$R$ is a Noetherian regular domain with precisely $n$ maximal 
  ideals, each maximal ideal of $R$ is of height 2, and the localizations of $R$ at its maximal
  ideals are the $n$ elements in $Q(D)$ that intersect irredundantly  to give $R$.   
  
  Let $R$ be a normal overring of $D$ such that each maximal ideal of $R$ has height 2. 
  Theorem~\ref{normal} asserts: (i)   $R$ is Noetherian if and only if $R$ is a flat overring of a 
  finitely generated $D$-subalgebra of $R$,  and (ii) if $R$ is Noetherian and local,  
  then $R$ is a spot over $D$. 
  
    Theorem~\ref{cor3.6}   asserts that the normal Noetherian overrings of $D$ with height 2 maximal
    ideals are precisely the rings $R$ for which there exists a nonsingular projective model $X$ over $D$
    and a subset $\U$ of the closed points of $X$ such that $R = \bigcap_{T \in \U}T$.
        

In   Section~\ref{sec8} we 
consider irredundant intersections of rings in $Q(D)$.  We prove in Theorem~\ref{prop4.2} that the representation of $D$ as the intersection of its first neighborhood rings is irredundant, and that if $\U$ is a proper subset of the set of all such rings, then the intersection of the rings in $\U$ is a flat extension of a regular finitely generated $D$-subalgebra of $F$ and hence $\U$ is an essential irredundant   representation of  the ring $\bigcap_{T \in \U}T$.  

When $D$ is Henselian we obtain our strongest result regarding irredundance.
  Let $\U$ be a set of pairwise incomparable rings in $Q(D)$.   Theorem~\ref{thm5.10} establishes that if $D$ 
is Henselian,  then the representation $\bigcap_{R \in \U}R$ is irredundant.  In Corollary~\ref{Henselian cor} we use this to show that for $D$ Henselian, every  Noetherian normal overring $R$  of $D$ for which each maximal ideal has height $2$ is    an irredundant intersection of the regular local rings in $Q(D)$ that are minimal with respect to containing $R$.

  

  \section{Preliminaries} 
  
   Our notation is as in Matsumura   \cite{Mat}.  Thus a local ring need not be Noetherian.
   We refer to Swanson and Huneke \cite{SH} for material on Rees valuation rings and 
   blowing up of ideals.
We refer to an extension ring $B$ of an integral domain $A$ as an {\it overring of} $A$ 
if $B$ is a subring of the quotient field of $A$.   A local ring $B$ is said to be a {\it spot} over $A$, 
if $B$ is a localization of a finitely generated $A$-algebra.

  We use the following definitions.
  
 \begin{definition} \label{2.3}
{\em Let $R$ be a  Noetherian  local integral domain and let $S$ be a local  overring of $R$.
\begin{enumerate}[(1)]
\item The {\it center} of $S$ on $R$ is the prime ideal $\m_S \cap R$ of $R$, where $\m_S$ denotes
the maximal ideal of $S$.
\item 
$S$ is said to {\it dominate} $R$ if the center of $S$ on $R$ is the maximal ideal of $R$,  that is, 
$\m_S \cap R = \m_R$, where $\m_R$ is the maximal ideal of $R$.

\item
If $\dim R \ge  2$,  a valuation overring $V$ of $R$ centered on $\m_R$  is said to be a  
{\it  prime divisor of the second kind } on $R$ if   
 the  field $V/\m_V$ has transcendence
degree $\dim R - 1$  over the field $R/\m_R$.\footnote{ See Zariski-Samuel    \cite[p.~95]{ZS2}. Valuation overrings of 
$R$ centered on height 1 primes are {\it prime divisors of the first kind}.   Prime divisors are necessarily DVRs.}
 
\item
$V$ is said to be a {\it  minimal valuation overring} of $R$ if $V$ is minimal with respect to set-theoretic
inclusion in the set of valuation overrings of $R$.  
\end{enumerate}}
  \end{definition}

\begin{remark} \label{rm1.2}  {\em 
Assume  notation as in Definition~\ref{2.3}. 
\begin{enumerate}[(1)]
\item  If $W$ is a valuation overring of $R$ and
the center $\m_W \cap R $ of $W$ on $R$  is a nonmaximal prime  ideal of $R$, then 
by composite  construction \cite[p.~43]{ZS2},  there exists a valuation overring $V$ of $R$ such that
$V \subset W$ and $\m_V \cap R = \m_R$.  Therefore every valuation overring of $R$ contains a valuation
overring of $R$ that is 
centered on the maximal ideal of $R$. 
\item   If $W$ is a valuation overring of $R$ that dominates $R$ and the 
field $W/\m_W$ is transcendental over $R/\m_R$,  then by composite construction,
 there exists a valuation overring $V$ of $R$ such that
$V \subset W$.  
\item 
 Every valuation overring of $R$ contains a minimal valuation overring of $R$.
 \item 
 Let $V$ be a valuation overring of $R$.  Then 
  $V $ is a minimal valuation overring of $R \iff    V$
dominates $R$ and   the field $V/\m_V$ is algebraic over the field $R/\m_R$.
\end{enumerate}}
\end{remark}
  
    Abhyankar in  Proposition~3 of \cite{Ab1}  characterizes prime divisors of a  regular local domain
    centered on the maximal ideal.
  The characterization is as follows.

 \begin{theorem} \label{1.20}
Let $R$ be a regular local domain  with  $\dim R = n \ge 2$ and let  $\m_R$ denote the maximal ideal
of $R$.  
Let $V$ be a prime divisor of   $R$ centered on $\m_R$.
 There exists a unique  finite sequence
\begin{equation} \label{1}
R ~ =  ~R_0  ~\subset ~ R_1 ~ \subset \cdots \subset ~ R_{h}  ~\subset ~ R_{h+1} ~= ~ V
\end{equation}
of regular local rings $R_j$, where  $\dim R_h \ge 2$ and    $R_{j+1}$ is the first local quadratic transform  of
$R_{j}$ along $V$   for each $j \in \{0, \ldots, h \}$,
and $\ord_{R_{h} }  = V$.\footnote{For the definition of  quadratic transforms, see for example \cite[pp.~569--577]{Ab8} and \cite[p.~367]{ZS2}. The powers of the maximal ideal of a regular local domain $S$ define a rank one
discrete valuation domain  denoted  $\ord_S$.   If $\dim S = d$,  then the residue field of $\ord_S$ is a pure transcendental
extension of the residue field of $S$ of transcendence degree $d-1$.}
\end{theorem}

It follows from Theorem~\ref{1.20}  that  the residue field $V/\m_V$ of $V$ is a pure
transcendental extension of
the field $R_h/\m_{R_h}$ of transcendence degree   one less than $\dim R_{h} $.
 Therefore the residue field of $V$
is ruled as an extension field of the residue field of $R$.\footnote{A field extension
 $F \subset L$ is said to be {\it ruled}
if $L$ is a simple transcendental extension of  a subfield $K$ such that  $F \subset K$.}
If  $\dim R \ge 2$ a prime divisor on $R$ is not a minimal valuation overring of $R$.

The association of the prime divisor $V$ with the regular local ring $R_h$ in Equation~\ref{1}, and the
uniqueness of the sequence in Equation~\ref{1}
establishes a one-to-one correspondence between the prime divisors $V$ dominating the
regular local ring $R$ and the regular local rings $S$ of dimension at least 2
 that  dominate $R$ and are
obtained from $R$  by a finite sequence of local quadratic transforms    as in Equation~\ref{1}.
The regular local rings $R_j$ with $j \le h$ displayed in Equation~\ref{1} are called 
the {\it infinitely near points } to $R$
along $V$.  In general, a regular local ring $S$ of dimension at least 2  is
called an {\it infinitely near point } to $R$ if
there exists a sequence
$$
R ~ =  ~R_0  ~\subset ~ R_1 ~ \subset \cdots ~ \subset ~ R_{h}~ = ~ S,  \quad h ~\ge ~ 0
$$
of regular local rings $R_j$ of dimension at least 2,   where   $R_{j+1}$ is the first
local quadratic transform  of
$R_{j}$   for each $j$ with $0 \le j \le h-1$
\cite[Definition~1.6]{L}.

\section{The quadratic tree of $D$}

Let $D$ be a $2$-dimensional regular local ring. 
 The Zariski-Abhyankar Factorization Theorem \cite[Theorem~3]{Ab1} implies that 
 every 2-dimensional regular local ring $R$ that birationally
dominates
 $D$ is an infinitely near point
to $D$. Because we will often be treating such rings  as points in what follows, we follow Lipman \cite{L} and  denote the infinitely near points to $D$ with Greek letters.  
    We record  in Theorem~\ref{1.1} implications of \cite[Theorem~3 and Lemma~12]{Ab1}.

\begin{theorem} \label{1.1} Let $D$ be a $2$-dimensional regular local ring, and let $\alpha$ be a $2$-dimensional regular local ring that birationally dominates $D$.   
\begin{enumerate}[(1)]
\item If $D \ne \alpha$, then $\m_D $ extends to a proper principal ideal of $\alpha$.
Therefore $\alpha$ dominates a unique local quadratic transform $\alpha_1$ of $D$.
\item There exists for some positive integer $\nu$ a sequence
$$
D~= ~ \alpha_0 ~\subset ~\alpha_1 ~\subset ~ \cdots ~ \subset \alpha_{\nu} ~= ~ \alpha,
$$
where $\alpha_i$ is a local quadratic transform  of $\alpha_{i-1}$ for each $i \in \{1, \ldots, \nu\}$.
The rings $\alpha_i$ are precisely the regular local domains  that are subrings of $\alpha$ and contain
$D$.
\item
If $V $ is a minimal valuation ring of $D$, then $V$ is the union of the infinite quadratic sequence of $D$ along $V$.

\end{enumerate}
\end{theorem}

 \begin{definition} \label{quad tree} {\em 
 Let $D$ be a regular local ring of dimension 2  and let $F$ denote the quotient field 
 of $D$.
 \begin{enumerate}[(1)]
 \item The {\it quadratic tree} 
   $Q(D)$ of $D$ is the partially ordered set (ordered by inclusion) defined as   the set
 of all iterated quadratic transforms of $D$.   Theorem~\ref{1.1}  implies that 
  $Q(D)$ is the set of all 2-dimensional 
 regular local rings that birationally dominate $D$.\footnote{This is the notation used by Abhyankar
 in  several  papers such as \cite{Ab14}.}
 
 \item 
 $Q(D)$ is the disjoint union of sets $Q_j(D)$ for $j \ge 0$, where $Q_j(D)$ denotes
 that set of all 2-dimensional regular local rings that 
 are  obtained  by making  precisely  $j$ quadratic transforms starting at $D$,  where $Q_0(D) = \{D\}$,  
  see  \cite{Ab1} and   \cite{Ab2}. We refer to the elements in $Q_j(D)$ as 
  {\it infinitely near points at level} $j$ to $D$.

  \end{enumerate}}
  \end{definition}
  
  \begin{remark} \label{rem3.3}
  
 {\em 
  \begin{enumerate}[(1)]
  \item[]
  \item   
  The elements in $Q_1(D)$ are in one-to-one correspondence with the points on a projective line over the 
  residue field $\kappa(D) = D/\m_D$ of $D$.  Therefore the  infinitely near points at level 1 to 
  $D$ are in one-to-one correspondence with the 
   irreducible  homogeneous polynomials in $\kappa(D)[x, y]$.  For each irreducible homogeneous 
   polynomial $f \in \kappa(D)[x,y]$ there exists an infinitely near point $\alpha_f \in Q_1(D)$.
   Irreducible homogeneous polynomials $f$ and $g$ in $\kappa(D)[x,y]$ are such that 
   $\alpha_f = \alpha_g$ if and only if  $f$ and $g$ are associates in $\kappa(D)[x, y]$.

  \item 
  Assume that $\kappa(D)$ is algebraically closed and that $D$ has a coefficient field, that is,
  there exists a subfield $k$ of $D$ that maps onto $ \kappa(D)$ under the natural surjection 
  $D \to D/\m_D$. Then  each  infinitely near points at level 1 to 
  $D$ is uniquely determined by a nonzero   homogeneous linear  polynomial in $k[x, y]$. 
  For $a, b \in k$ with  $a \ne 0$,  the polynomial   $ay  + bx $ is associated to 
  $D[\frac{y}{x}]_{(\frac{ay + bx}{x})}$.  If $a = 0$,  then $b \ne 0$ and we may assume $b = 1$
  and associate $y$ to the local quadratic transform
  $D[\frac{x}{y}]_{(\frac{x}{y})}$.

 

 \end{enumerate}
}
 \end{remark}


For future reference, we collect here notation we will use throughout the article. 

\begin{notation} \label{not1.3}
{\em We use the following notation. 
\begin{enumerate}[(1)]

\item $D$ is a $2$-dimensional regular local ring with quotient field $F$ and maximal ideal ${\ff m}_D = (x,y)D$. 

\item $Q(D)$ is the quadratic tree of $D$ as in Definition~\ref{quad tree}. 

\item For each $\alpha \in Q(D)$ and $j \geq 0$, $Q_j(\alpha)$ is the set of infinitely points at level $j$ to $\alpha$.  
\item    
 For each   subset $\U$ of $Q(D)$,  let $\O_\U  = \bigcap_{R \in \U}R$. 
 
 \item
 Let $\R(D)$ denote the set of rings of the form $\O_\U$ for some subset $U$ of $Q(D)$.

  \end{enumerate}}
  \end{notation}

\begin{remark}  \label{disc2.7}  {\em   The Noetherian rings in $\R(D)$ are all Krull domains.  Associated to a Krull domain 
$A$ is a unique set of DVRs,  the set $\E(A)$  of  essential valuation rings of $A$; 
 $\E(A) =  \{A_\p \}$,  where $\p$  varies 
over  the  
height 1 prime ideals  of $A$.  Two useful properties related to $\E(A)$   are:   
\begin{enumerate}[(1)]
\item
$A = \bigcap \{ A_\p~ |~ A_\p \in \E(A) \}$ and the intersection is irredundant.
\item
The set $\E(A)$ defines an essential representation of $A$.

\end{enumerate} }

\end{remark} 

One of our motivations for this article and \cite{HLO} is to examine the extent to which there  are similarities between the intersections of elements in $Q(D)$ with the representation of a Krull domain $A$  as an 
intersection of  its   essential valuation rings.

 \section{ Projective models over $ D$} 
  

Let $D$ with quotient field $F$  be as in Notation~\ref{not1.3}. 
In this section we   relate  the geometry of $Q(D)$  to nonsingular projective models over $D$.
We use  the following terminology as in 
Section~17,  Chapter VI of 
Zariski-Samuel \cite{ZS2}.   
 If   $A$ is a finitely generated 
$D$-subalgebra of
 $F$, the   {\it affine model} over $D$  associated to $A$  is the set of local rings 
$A_\p$, where  $\p$ varies over the set of prime ideals of $A$. 
A  {\it model}  $M$   over
$D$ is a  subset of  
the local overrings of  $D$ that has   the properties:
(i)  $M$ is 
a finite union of affine models over $D$,  and (ii) each valuation overring of $D$ 
 dominates at most one of the local rings in $M$.  This second condition is called {\it irredundance}.  
  A model $M$ over $D$ is said to 
be {\it complete} if each valuation overring of $D$ 
dominates a local ring in $M$.

 A model  $M$ is said to be {\it projective} over $D$ if  there exists a finite
set of nonzero elements  $a_0, a_1, \ldots, a_n$ in $ D$ such that 
$J=(a_0,\ldots,x_n)D$ is an ${\ff m}$-primary ideal of $D$  and $M$ is the union 
of the affine models defined  by the rings 
$A_i  = D[\frac{a_0}{a_i}, \frac{a_1}{a_i}, \ldots, \frac{a_n}{a_i}],  i = 0, 1, \ldots, n$.  
This projective model is the {\it blowup}\footnote{We are identifying  the projective scheme $\Proj D[Jt]$ 
with the {model} $\bigcup_{i=0}^n \{(A_i)_{\ff p} ~|~{\ff p} $ in $\Spec A_i\}$.}   
 $\Proj D[Jt]$ of the ideal $J$. 
In the language of schemes, the projective and affine   models  
we consider correspond   to the projective schemes over $\Spec D$ and the affine schemes
 over $\Spec D$  of finite type. 
 
 A basis for the Zariski topology on a model $M$ is given by the sets of the form $\{R \in M:x_1,\ldots,x_n \in R\}$, where $x_1,\ldots, x_n \in F$.  The closed points in this topology are the local rings in $M$ that are maximal with respect to set inclusion. If $M$ is a projective model over $D$, then the closed points are precisely the local rings in $M$  of dimension~$2$.

 A model $M$ over $D$  is {\it normal} if every local ring in $M$ is a normal domain. 
 The {\it normalization} of a projective model $M$ over $D$  is the projective model over $D$ obtained by normalizing each affine component in $M$; i.e., if $M$ is the union of the affine models defined by the rings $A_0,\ldots,A_n$, then the normalization of $M$ is the union of the affine models  defined by the normalization of the  $A_i$.  If $M = \Proj D[It]$ is a projective model over $D$ for an ideal $I$ of $D$, then the normalization of $M$ is $\Proj D[Jt]$, where $J$ is the integral closure of the ideal $I$.\footnote{Here we are using that the powers of $J$ are also integrally closed, and $D[Jt]$ is a normal domain.} Thus a projective model $M = \Proj D[Jt]$ is normal if and only if $J$ is a complete ideal\footnote{Also called integrally closed
ideals.}.

Classical results proved by Zariski on the structure of complete ideals of a 2-dimensional regular local ring $D$  simplify the structure of projective models birational 
over $D$.  Complete ideals of $D$ are closed with respect to ideal multiplication,
and there is a marvelous unique factorization theorem:  every nonzero  complete ideal can be written 
uniquely as a  finite product of simple complete ideals,   cf.  \cite[Appendix 5]{ZS2}  
or \cite[Chapter 14]{SH}. 

We use the following terminology.  

\begin{definition}  \label{def4.1} {\em
Assume Notation~\ref{not1.3} and let $I$ be an $\m_D$-primary ideal. 
\begin{enumerate}[(1)]
\item
The {\it base points} of $I$  
are the points $\alpha \in Q(D)$ for which the transform of $I$ in $\alpha$ is a proper ideal of $\alpha$. 
Let  $\B(I)$   denote  the set of base points of $I$.   Then 
 $\B(I)$ is a finite  subset of $Q(D)$    \cite{L}.  A base point $\alpha$ of $J$ is called a {\it maximal} or {\it terminal} base 
 point of $J$ if $\alpha$ is a maximal element of the partially ordered set $\B(J)$, cf.  \cite[Remark~2.9]{HKT2}.
\item
The set $\Rees I$ of {\it Rees valuation rings}  of $I$ is the  smallest set $\{V_1,\ldots,V_n\}$ of  valuation 
overrings of $D$
such that for each $k>0$, the integral closure of $I^k$ is $I^kV_1 \cap \cdots \cap I^kV_n \cap D$.  
The set with this property is unique, and each $V \in \Rees I$ is the order valuation ring 
$\ord_\alpha$ of a unique point $\alpha \in Q(D)$. 
\item
Let $J$ be a simple complete $\m_D$-primary ideal. Then $\Rees J = \{ \ord_\alpha \}$ for a 
unique point $\alpha \in Q(D)$, cf. \cite[Prop.~14.4.8]{SH}.     As in Theorem~\ref{1.20},  there exists a 
unique chain $D = \alpha_0 \subset \alpha_1 \subset \cdots \subset \alpha_n = \alpha$ of infinitely
near points from $D$ to $\alpha$. 
Then $\B(J) = \{\alpha_0, \ldots, \alpha_n\}$ is the set of base points of $J$.  There exists a descending
sequence $\m_D = J_0 \supset J_1 \supset \cdots \supset J_n = J$ of simple complete ideals 
of $D$,  where $\Rees J_i = \{\ord_{\alpha_i} \}$ for each $i$.  The {\it saturation of }  $J$ is 
the ideal $L = \prod_{i=0}^n J_i$.  
\item
To define the saturation of an $\m_D$-primary ideal $I$,   let $J$ be the integral closure of $I$. 
For each simple complete factor $J_i$ of $J$, let $L_i$ denote the saturation of $J_i$. 
The {\it saturation} $L$ of $I$ is the product of the ideals $L_i$ as we vary over all the distinct 
simple complete factors of $J$.\footnote{Two ideals with the same simple complete factors define the same blowup. Thus
  in defining the saturation of an ideal, it does not matter if a given complete simple factor occurs more than once.}
 
\item
If $L$ is the saturation of $I$, then $L$ is also the saturation of $L$,  and we say that 
$L$ is a {\it saturated ideal}.

\end{enumerate}}
\end{definition} 

  We summarize in Remark~\ref{facts2.3} properties of saturated ideals that follow from the definition.

 \begin{remark}    \label{facts2.3}
{\em  Assume Notation~\ref{not1.3},   and let $J$ be a complete 
  $\m_D$-primary ideal. 
  \begin{enumerate}[(1)]
\item
 Assume  $V  \in \Rees J$.   Then $V$ is the order valuation ring of  
  $\alpha_n \in Q_n(D)$,  for some integer $n \ge 0$.
 Let $D \subset \alpha_1 \subset \alpha_2 \subset \cdots \subset \alpha_n$ be the unique chain of regular local rings 
 from $D$ to $\alpha_n$.   If $J$ is saturated, then  the order valuation rings for $D, \alpha_1, \ldots, \alpha_{n-1}$ are 
  in the set $\Rees J$
    \item
    $\Rees J \subseteq \{\ord_\alpha ~|~ \alpha \in \B(J)\}$.
 \item 
  $J$ is saturated  $\iff   \Rees J  =  \{\ord_\alpha ~|~ \alpha \in \B(J) \} $.  
  \end{enumerate} }
 \end{remark}

In  \cite[Definition 5.11]{HJLS}  the following equivalent formulation 
to Definition~\ref{def4.1} 
of a 
saturated ideal is given. 

\begin{remark}  \label{rem4.2}
{\em  A complete   $\m_D$-primary ideal $J$ 
is  saturated  if for each 
simple complete ideal $I$ with $J \subseteq I$ and $I = IV \cap D$ for some $V \in \Rees J$, the ideal $I$ is a factor of $J$.   }
\end{remark}

\begin{pf}
 The equivalence follows because if $V \in \Rees J$ and 
$V = \ord_\alpha$,  then $V$ dominates $\alpha$ and therefore $V$ dominates each of the infinitely near 
points in the chain from $D$ to $\alpha$. Let $I_n$ be the simple complete ideal corresponding to $\alpha$.
Then  $V \in \Rees J$ implies that 
 $I_n$ is a factor of $J$ by the unique factorization theorem of Zariski \cite[Theorem~14.4.9]{SH}. 
Moreover,  the simple complete ideals corresponding to points in the chain from 
$D$ to $\alpha$  are contracted from $V$. The condition in Remark~\ref{rem4.2} implies  that all these 
simple complete ideals are also factors of $J$. Hence $J$ is saturated.  

Conversely,  if $J$ is saturated, then $\Rees J = \B(J)$,  and the condition in Remark~\ref{rem4.2} holds. 
\qed \end{pf}

  Saturation has an important geometric interpretation.  
 A model $M$ over $D$ 
 is {\it nonsingular} if every  ring in $M$ is a regular local ring.  
 We record in Proposition~\ref{sat = ns} a result given in \cite[Proposition 5.12]{HJLS}.
 
 \begin{proposition} \label{sat = ns}   A normal projective model $M = \Proj D[Jt]$ over $D$, 
 where $J$ is a complete ${\ff m}_D$-primary ideal, is nonsingular if and only if $J$ is saturated. 
 \end{proposition}

Facts~\ref{disc3.6}  records known properties of a nonsingular projective model
$X$   over $\Spec D$. In our description of $X$ it is useful that $X = \Proj  D[Jt]$, 
where $J$ is a saturated complete ${\ff m}_D$-primary ideal. We are interested in 
relating  properties of $X$ to the quadratic tree  of $D$. The finite set $\B(J)$ of 
base points of $J$ plays an important role in this connection.  

\begin{facts} \label{disc3.6}
{\em  Assume Notation~\ref{not1.3}.  Let $X$ be a nonsingular  projective model 
 over $\Spec D$ and let $J$ be a saturated complete 
 $\m_D$-primary ideal such that $X = \Proj D[Jt]$. 
 The closed points of $X = \Proj D[Jt]$  are a subset of $Q(D)$ of a special form that is related to the set $\Rees J$ as follows:
 \begin{enumerate}[(1)]
 \item 
 All but finitely many of the points in $Q_1(D)$ are in $X$. 
 \item
 The following are equivalent:
 \begin{enumerate}[(a)]
\item
All the  points in $Q_1(D)$ are in $X$.
 \item   
 $Q_1(D)$ is the set of closed points of   $X$. 
 \item
   $ \B(J)  =  \{D\}$.
   \item
   $  \{ \ord_D \}   = \Rees J $.
   \item
   $  J  = {\ff m}^k_D$ for some $k>0$.
   \end{enumerate}
 
 \item
 Let $\alpha \in \B(J)$.  Then:
 \begin{enumerate}[(a)]
 \item
 Each of the points  of $Q(D)$ in the chain from $D$ to $\alpha$ is in $\B(J)$.
 \item 
  $\alpha \not\in X$.
  \item
  $\ord_\alpha \in \Rees J$.
  \item
  All but finitely many of 
 the points of $Q_1(\alpha)$ are in $X$.
 \end{enumerate}
 \item 
 $\alpha \in \B(J)$ is a terminal base point of $J$  as in Definition~\ref{def4.1}  $\iff Q_1(\alpha) \subset X$.
 \item 
 The finite set of terminal base points of $J$
 uniquely determines the nonsingular projective model $X = \Proj D[Jt]$.  
 \item
 There exists an integer $s \ge 0$ such that the terminal base points of $J$ are 
 contained in  $Q_0(D) \cup Q_1(D) \cup \cdots \cup Q_s(D)$.  
 \item
 If the terminal base points of $J$ are 
  contained in  $Q_0(D) \cup   \cdots \cup Q_s(D)$, 
  then   the closed points of $X$ are contained in $Q_1(D) \cup Q_2(D) \cup \cdots \cup Q_{s+1}(D)$.
 \end{enumerate} }
\end{facts}


Facts~\ref{disc3.6}  implies the following characterization of the closed points of a nonsingular projective
model $X$  over $\Spec D$  in terms of elements of the quadratic tree $Q(D)$, and the 
base points $\B(J)$ of $J$.

\begin{theorem} \label{cor5.2}
Let $X = \Proj D[Jt]$ be a nonsingular projective model over $D$  as in Facts~\ref{disc3.6}.
The set $\U$ of closed points of $X$  has the following form:  
Either
 \begin{enumerate}[(1)]
 \item
 $\U = Q_1(D)$ in which case  $\B(J) = \{ D \}$ and  $X = \Proj D[xt, yt]$, or
 \item
   $   \B(J)  \setminus \{D\} = \{\alpha_1,\ldots,\alpha_n \} $,  and  
  $$\U = \left(Q_1(D) \cup Q_1(\alpha_1) \cup \cdots \cup Q_1(\alpha_n)\right) \setminus \{\alpha_1,\ldots,\alpha_n\},$$
 in which case $X = \Proj D[Jt]$ for a saturated complete    $\m_D$-primary ideal  $J$  such that  $\Rees J$ is the set of order valuation rings of the rings in $\B(J)$.   
 \end{enumerate}
 Conversely, if $\S  =   \{\alpha_1,\ldots,\alpha_n \} $ is a   finite subset of $Q(D)  \setminus \{D\}$ 
 having the property that $\alpha \in \S$ implies each point in the chain for $D$ to $\alpha$ is 
 in the set $\S \cup \{ D\}$,  then there exists a saturated complete ideal $J$ such that 
  $$\U = \left(Q_1(D) \cup Q_1(\alpha_1) \cup \cdots \cup Q_1(\alpha_n)\right) \setminus \{\alpha_1,\ldots,\alpha_n\}$$
  is the set of closed points of $X = \Proj D[Jt]$.
\end{theorem}

\begin{pf} Apply Facts~\ref{disc3.6}.   
\qed \end{pf}

\section{Desingularization of projective models}

As Proposition~\ref{sat = ns} suggests, saturation is the algebraic analogue of desingularization. We formalize this connection in Theorem~\ref{thm2.6}. 
We recall  the desingularization of a projective model, as defined in \cite[p.~199]{L1969}. 


\begin{definition} \label{def2.5} {\em  Let  $M$ and $N$ be  models over $ D$.   Then $N$ {\it dominates} $M$ if  each valuation overring   $V$ of $D$ centered on a ring in  $N$ dominates the center of $V$ on $M$; equivalently, 
 each local ring in $N$ dominates a local ring in $M$.\footnote{Viewing $N$ and $M$ as projective schemes over $\Spec D$, this implies  there is a  birational morphism $N \rightarrow M$.}

 Let  $R$ be a  Noetherian overring of $D$.  
 Let $M$ be a projective model over $R$.  
 \begin{enumerate}[(1)]
 \item
 A {\it desingularization}\footnote{Our definition differs  from Lipman's but is equivalent in our context.  Following \cite[p.~199]{L1969}, 
    a desingularization of a projective model $M$ of $D$ is 
     a proper birational map of surfaces, $N \rightarrow M$, such that $N$ is nonsingular. In this case, $N \rightarrow \Spec D$ is also a proper birational map of surfaces. This fact, along with the assumptions that $N$ is nonsingular and $D$ is a Noetherian ring, implies $N$ is projective  \cite[Corollary 27.2]{L1969}. Thus  our definition is equivalent to Lipman's in our setting.}  
 of $M$ is a nonsingular projective model $N$ over $D$ that dominates $M$.
 \item
 A desingularization $N$ of $M$ is a {\it minimal desingularization}\footnote{ If there exists a desingularization  $N$ of $M$, then there exists a unique minimal desingularization of $M$ \cite[Corollary 27.3]{L1969}. }  if every desingularization of $M$ 
 dominates $N$.
 \end{enumerate}}
 
\end{definition}

The Zariski theory of complete ideals  along with the Zariski-Abhyankar factorization theorem 
yields the following result. 

\begin{theorem}  \label{thm2.6}
Let $J$ be  a complete ideal of $D$,  let $M = \Proj D[Jt]$ and let   
 $L$ denote the saturation of $J$. 
Then $N = \Proj D[Lt] $ is a minimal desingularization of $M$.
The converse also holds;  if $L'$ is a complete ideal such that  
$\Proj D[L't]$ is a minimal desingularization of $M$,  
then 
$L$ and $L'$ have the same simple complete factors and 
 $\Proj D[L't] = \Proj D[Lt]$.
\end{theorem}


\begin{pf}  {The model $N$  dominates $M$. Also, 
$N = \Proj D[Lt]$ is a 
 nonsingular model over $ D$ by Proposition~\ref{sat = ns}.}  Hence $ \Proj D[Lt]$ is a 
desingularization of $\Proj D[Jt]$.

Since     $\ord_{\alpha} \in \Rees J$ for each terminal base point  $\alpha \in \B(J)$,
the set $\B(J)$ of base points of $J$ is the same as the set $\B(L)$ of base points of $L$. 
 Theorem~\ref{cor5.2}  implies that $\Proj D[Lt]$ is the unique nonsingular projective model
over $D$ having  the set $\B(L) = \B(J)$ as base points.

Let $Y$ be a nonsingular projective model over $D$ that dominates  $\Proj D[Jt]$.  
 There exists a complete saturated ideal $I$ of $D$ such that $ Y =   \Proj D[It]$. 
 
 If $V \in \Rees J$, then $V \in \Proj D[Jt]$.  
 Since $\Proj D[It]$ dominates $\Proj D[Jt]$
  {and the only local ring birationally dominating a  
   valuation ring is the valuation ring itself}, it follows that $V \in \Proj D[It]$ for  each $V \in \Rees J$. 
 Therefore  
  the set $\B(I)$ of base points of $I$ 
 contains $\B(J) = \B(L)$. 
 Proposition~\ref{sat = ns} implies 
  that each of the simple complete factors of $L$ is a 
 factor of $I$.   Therefore $Y = \Proj D[It]$ dominates $\Proj D[Lt]$,  and the domination map
 $Y \to \Proj D[Jt]$ factors through $Y \to \Proj D[Lt]$.
 
 For the converse,  if $L'$ is a complete ideal such that  $\Proj D[L't]$ is a minimal desingularization of $M$, 
 then $\B(L') = \B(J) = \B(L)$,  and the complete ideals $L$ and $L'$ have the same complete 
 simple factors. 
  Therefore $\Proj D[Lt] = \Proj D[L't]$ and $L'$ is also the saturation of $J$.
\qed \end{pf}

While Theorem~\ref{thm2.6} characterizes the desingularization of a normal projective model $M = \Proj D[Jt]$ in 
terms of the saturation of the ideal $J$,  Theorem~\ref{order-theoretic} 
 characterizes the desingularization of $M$ strictly in terms of order-theoretic properties of $Q(D)$.  

\begin{theorem}  \label{order-theoretic}
Let $M$ be a normal projective model over $D$. The closed points of the minimal desingularization of $M$ are the points in $Q(D)$ that are minimal with respect to dominating a closed point in $M$.  
\end{theorem}

\begin{pf}  We may assume $M$ has singularities since otherwise the theorem is clear. By Theorem~\ref{thm2.6} there is a unique minimal desingularization $N$  of $M$.  Since $M$ is a  normal surface that can be desingularized, $M$ has finitely many singularities \cite[Theorem, p.~151]{L1969}. By \cite[Propositions~1.2 and~8.1]{L1969} and  \cite[B, p.~155]{Lip2},  there is a sequence $M_n \rightarrow \cdots \rightarrow M_1 \rightarrow M_0 = M$ of normal projective models over $D$ such that $M_n$ is nonsingular and for each $i$, $M_{i+1}$ is obtained from $M_i$ by blowing up the finitely many singular points of $M_i$.  
Let $\U$ be the set of points in $Q(D)$ that are minimal with respect to dominating a closed point in $M$.  
We claim that $\U$ is the set of closed points of $M_n$.  

Let $\alpha \in \U$, and let $R$ be the center of $\alpha$ in $M$.  If $\alpha \ne R$, then $R$ is a singular point in $M$. Since $\alpha$ dominates $R$ and $R$ is a singular point in $M$, $\alpha$ dominates a point in $M_1$ \cite[($\star$), p.~203]{L1969}.  If $\alpha \not \in M_1$, then  $\alpha$ dominates a singular point in $M_1$ since $\alpha \in \U$. 
Continuing in this manner, we obtain either that $\alpha \in M_i$ for some $i$ or $\alpha \not \in M_n$ and $\alpha$ dominates a point in $M_n$. The latter property is contrary to the fact that $M_n$ is nonsingular and $\alpha$ is minimal among points in $Q(D)$ dominating $R$.   Thus $\alpha \in M_i$ for some $i$, and since $\alpha$ is a nonsingular point in $M_i$ and in the sequence $M_n \rightarrow \cdots \rightarrow M_1$ we have only blown up singular points, we have $\alpha \in M_n$.  This shows that every point in $\U$ is a closed point in $M_n$.  

Conversely, let $\alpha$ be a closed point in $M_n$, and let $R$ be the center of $\alpha$ in  $M$.   Let $\beta \in Q(D)$ such that $\beta \subseteq \alpha$ and $\beta$ is  minimal with respect to dominating $R$.  If $\beta \in M$, then $\beta$ is a nonsingular point in $M$. By the construction of the $M_i$, it follows then that $\beta \in M_n$, so that $\beta = \alpha$. Otherwise, $\beta \not \in M$ and so $R$ is a singular point in $M$. Thus $\beta$ dominates a point $R_1$ in $M_1$ \cite[($\star$), p.~203]{L1969}.  Continuing in this manner and using the fact that $\alpha \in M_n$ dominates $\beta$, we obtain that $\beta \in M_i$ for some $i$ and hence $\beta = \alpha$. Therefore, $\alpha$ is minimal with respect to dominating its center in $M$, which proves that every closed point in $M_n$ is in $\U$.  

Now consider the minimal desingularization $N$ of $M$.  Each closed point in $N$ dominates a point in $\U$, so $N$ dominates $M_n$. Since $N$ is a minimal desingularization of $M$, we conclude that $N = M_n$, which complete the proof.
\qed \end{pf}

\begin{remark} \label{join5.4} {\em
 Let $M$ and $M'$  be projective models over $D$.   We refer to \cite[p.~120]{ZS2} 
 for the definition of the join $M''$ of $M$ and $M'$.  $M''$ is a projective model over $D$ that dominates
 both $M$ and $M'$.  If   $I$ and $I'$ are nonzero ideals of $D$
 such that  $M = \Proj D[It]$  and $M' = \Proj D[I't]$,  then $M'' = \Proj D[II't]$.  
 
 Assume that $M$ and $M'$ are nonsingular projective models over $D$.  There
 exist saturated complete $\m_D$-primary ideals $I$ and $I'$ such that 
 $M = \Proj D[It]$ and $M' = \Proj D[I't]$.  The ideal $II'$ is also complete and saturated:  
 $\B(II') = \B(I) \cup \B(I')$ and $\Rees(II') = \Rees I \cup \Rees I'$,  cf.  \cite[Prop.~10.4.8]{SH}.  
 Hence the join $M'' = \Proj D[II't]$ is also nonsingular. The set $\U''$ of closed points of $M''$
 has the following description in terms of the sets $\U$ and $\U'$ 
  of closed points of $M$ and $M'$ and the base points $\B(I)$ and $\B(I')$:  $\U'' = (\U \cup  \U' )
  \setminus (\B(I) \cup \B(I'))$.  }
\end{remark}

By Fact~\ref{disc3.6}.5, a finite set $\alpha_1,\ldots,\alpha_n$ of incomparable rings in $Q(D)\setminus \{D\}$ uniquely determines a nonsingular projective model $M$ over $D$ that has  the $\alpha_i$ as
 precisely the terminal base points of $M$. Since the points $\alpha_i$ are base points of $M$,
 they are not points of $M$.   On the other hand, Theorem~\ref{minimal2}  describes nonsingular 
 projective  models $M$  in 
 terms of finitely many points $\alpha \in M$. 
 
 As an application of saturation and desingularization,  
  Theorem~\ref{minimal2}   records properties of nonsingular models over $D$ in 
  terms of finite subsets of incomparable points of $Q(D) \setminus \{D\}$.

 
\begin{theorem} \label{minimal2}   Assume Notation~\ref{not1.3}.
\begin{enumerate}[(1)]
\item
Let $\alpha \in Q(D) \setminus \{D\}$. Then:
\begin{enumerate}[(1)]
\item
There exists a unique nonsingular projective model $M = \Proj D[Lt]$  such that $\alpha \in M$ and 
every nonsingular projective model over $D$ that contains $\alpha$ dominates $M$.
\item
The set $\U$ of closed points of $M$ has the property that  $\U \setminus \{\alpha\}$  is the set of 
points of $Q(D)$ minimal with respect to being incomparable to $\alpha$.
\end{enumerate}

\item
Let $\S  =  \{\alpha_1,\ldots,\alpha_n\}$ be  a finite   set of    incomparable points in $Q(D) \setminus \{D\}$. Then:
\begin{enumerate}[(1)] 
\item
There is a unique nonsingular projective model $M$ over $D$ such that $\alpha_1,\ldots,$ $\alpha_n \in M$ and every nonsingular projective model over $D$ containing $\alpha_1,\ldots,\alpha_n$ dominates $M$.  

\item
The set $\U$ of closed points of $M$ has the property that  $\U \setminus \S$  is the set of 
points of $Q(D)$ minimal with respect to being incomparable to every $\alpha \in \S$.
\end{enumerate}   
\end{enumerate}
 
\end{theorem}

\begin{pf}
For item~1,  let $\gamma \in Q(D)$ be the unique point 
such that $\alpha \in Q_1(\gamma)$.  Let $J$ be the 
simple complete ideal of $D$ having $\ord_{\gamma}$ as its Rees valuation ring,  and let $L$
be the saturation of $J$.  Let $M = \Proj D[Lt]$.  
By Fact~\ref{disc3.6}.5 and Theorem~\ref{thm2.6},  $M = \Proj D[Lt]$ satisfies item1.a.

Notice that $\gamma$ is the unique terminal base point of $M$,
and $\gamma \in Q_d(D)$ for some integer $d \ge 0$.  If $d = 0$,
then $\gamma = D$ and   $M = \Proj D[xt, yt]$.   It is clear in 
this case that $\U = Q_1(\gamma) \setminus \{\alpha\}$ is the set of points in $Q(D)$ minimal with respect to being 
incomparable to $\alpha$.  

Assume that $d \ge 1$ and let $D = \gamma_0 \subset \gamma_1 \subset \cdots \subset \gamma_d  = \gamma$
be the unique chain of infinitely near points from $D$ to $\gamma$.
As in Definition~\ref{def4.1}.3,  $\B(L) = \{\gamma_0,  \ldots \gamma_d \}$,
and the set $\U$ of closed points of $M = \Proj D[Lt]$ is 
 $$
 \U = \left(Q_1(D) \cup Q_1(\gamma_1) \cup \cdots \cup Q_1(\gamma_n)\right) \setminus \{\gamma_1,\ldots,\gamma_d\}.
 $$
 It follows also in this case that  $\U  \setminus \{\alpha\}$ is the set of points in $Q(D)$ minimal with respect to being 
incomparable to $\alpha$.  This verifies item 1.

For item~2,   for each  $\alpha_i \in \S$, there exists a point $\gamma_i$ such that 
$\alpha_i \in Q_1(\gamma_i)$.\footnote{For distinct $\alpha_i$ and $\alpha_j$, it may happen 
that $\gamma_i = \gamma_j$.}
By item~1,  for each $\alpha_i$, there exists a   nonsingular projective model $M_i = \Proj D[L_it]$
such that   $M_i$ has a unique terminal base point $\gamma_i$ 
and $\alpha_i \in Q_1(\gamma_i) \subseteq \U_i$,  where $\U_i$ is the set of 
closed points of $M_i$.

 Let $L$ be the product of the saturated complete ideals $L_i$.  Then $L$ is a
saturated complete ideal and $M = \Proj D[Lt]$ is the join of the models $M_i$. 
Theorem~\ref{thm2.6} and Fact~\ref{disc3.6}.5  imply that $M$ satisfies item~2.a. 

The set $\U$ of closed points of $M$ is
 $\bigcup_{i= 1}^n\U_i   \setminus \B(L)$,  where $\B(L)   = \bigcup_{i=1}^n\B(L_i)$ is the set of base points of $L$. 
Then   $\S \subset \U$ and $U \setminus \S$ is the set of 
points in $Q(D)$ minimal with respect to being incomparable  to every $\alpha \in \S$.
\qed \end{pf}

\section{Intersections of closed points in affine models}

In this section   we consider  the intersection of closed points in affine components  of nonsingular projective models over $D$. We  save the more subtle non-affine case for Section~\ref{sec7}.  

We use the following terminology  in  Theorem~\ref{prop3.7}. 

\begin{definition}   \label{def2.4}
{\em Let $\U$ be a subset  of $Q(D)$. A ring $\alpha \in \U$ is {\it essential} for $\U$ if  $\alpha$ is a localization of $\O_\U = \bigcap_{\alpha \in \U}\alpha$. 
We say  $\U$ 
 defines an  {\it  essential representation} of $\O_\U$ if each 
$\alpha \in \U$ is essential for $\U$. }
\end{definition}

\begin{remark} \label{rm2.5}
{\em Let $\U' \subset \U$ be subsets of $Q(D)$. 
If $\U$ defines an essential representation 
of $\O_\U$,    then $\U'$ defines an essential representation of $\O_{\U'}$.}
\end{remark}

\begin{pf}  Let $B = \O_\U$ and $C = \O_{\U'}$.  Then  $\U' \subset \U$ implies $B \subseteq C$.
Let $\alpha \in \U'$ and let $\m_\alpha$ denote the maximal ideal of $\alpha$.  
Since $\U$ defines an essential representation of $B$,  we have 
$ \alpha   = B_{\m_\alpha \cap B}$.  Then 
$ \alpha   = B_{\m_\alpha \cap B} ~ \subseteq C_{\m_\alpha \cap C} ~ \subseteq ~ \alpha$  implies 
$C_{\m_\alpha \cap C} ~ = ~ \alpha$.  Therefore $\U'$ defines an essential representation of $C$. 
\qed \end{pf}

\begin{theorem} \label{prop3.7}
Assume that $X = \Proj D[Jt]$   is a nonsingular projective model over $ D$. 
  Let $b \in J $ be part of a minimal set of generators of $J$ and 
 let $A = D[J/b]$ be the  associated affine component of $\Proj D[Jt]$.  Then   $A$ is a regular Noetherian domain.  Let $\U \subset Q(D)$  denote the set of  closed points of $X$ that contain $A$, and 
 let $B = \O_\U$.  
 \begin{enumerate}[(1)]
 \item 
 Each $\alpha  \in \U$ 
 is a localization of $A$ at a height 2 maximal ideal.
 \item 
 $A_\m \in \U$  and $A_\m = B_{(\m A_\m \cap B)}$ 
  for each height 2 maximal ideal $\m$ of $A$. 
  \item
 $\U$ defines an essential representation of $\O_\U$.
 \item
 If $\q$ is a maximal ideal of $B$, then $\hgt \q = 2$ and $B_\q = A_{\q \cap A}$.
  
 \item     
  $ \O_\U = B $ is a flat overring of  $A$.
\end{enumerate}
 \end{theorem}
\begin{pf}
Since $A$ is an affine component of $\Proj D[Jt]$,  $A$ is a regular Noetherian domain, and each $\alpha \in \U$ 
is a localization of $A$ at a height 2 maximal ideal. Let $\m$ be a maximal ideal of $A$. 
Then $\hgt \m = 2    \iff \m \cap D = \m_D \iff A_\m \in \U$.   If $A_\m = \alpha \in \U$,
then $\alpha$  is a localization of $B  = \O_\U$.   
This proves items~1,  2 and 3.   

Let $\q$ be a maximal ideal of $B$ and let $\p = \q \cap A$.  Let $\m$ be a maximal ideal
of $A$ with $\p \subseteq \m$.  If $\hgt \m = 2$, then $A_\m = B_\n$,  where
$\n = \m A_\m \cap B$.  It follows that $A_\p = B_\q$.  Since $\q$ is a maximal 
ideal of $B$, $\q = \n$ and $\m = \p$ in this case. 

If $\hgt \m = 1$,  then $\p = \m$ is a maximal ideal of $A$.  Since $A$ is
Noetherian and $A_\p$ is a DVR,  it follows that $\p$ is an invertible 
ideal of $A$.  Let $\p^{-1}$ denote the inverse of $A$.
Since $ B ~ =  ~  \bigcap\{A_\m ~|~ \m \in \Spec A,  ~\hgt \m = 2  ~ \}$,  we
have $\p^{-1} \subset B$ and $\p B = B$.  
Therefore this case does not occur.  This proves item~4.  
Theorem 2 of \cite{Ric}  now  proves item~5.
\qed \end{pf} 

\begin{discussion}\label{dis3.3} {\em We illustrate Theorem~\ref{prop3.7} in a special case. 
Assume Notation~\ref{not1.3}.  Then  $\Proj D[xt, yt]$  is the nonsingular projective model over $ D$
obtained by blowing up $\m_D$.  The closed points of  $\Proj D[xt, yt]$  are precisely the points of $Q_1(D)$.  
Let $\beta = D[x/y]_{(y, x/y)D[x/y]}$  denote the point in $Q_1(D)$ in
the $y$-direction.  Then $D[y/x]$ is the affine component  of  $\Proj D[xt, yt]$ that omits the point $\beta$.
The localizations of $D[y/x]$ at its height 2 maximal ideals describe the set $\U = Q_1(D) \setminus \{\beta\}$.
Let $A = \O_\U$.  Then $A$ is a localization of $D[y/x]$ at the multiplicatively closed set generated by 
principal generators of the height one maximal ideals of $D[y/x]$. There are infinitely many 
height one maximal ideals of $D[y/x]$.  For each integer $n \ge 2$,  let $V_n = D_{(y^n - x)D}$.   The 
center of $V_n$ on $D[y/x]$ is a height one maximal ideal of $D[y/x]$. 
 Notice that 
$xD[y/x]  = \m_D D[y/x]$ is a principal height one
prime ideal  that is contained in every height 2 maximal ideal of $D[y/x]$. It follows that $xA$ is the 
Jacobson radical of $A$  and $A = (1 - xD[y/x])^{-1}D[y/x]$.  Also $D[y/x]/xD[y/x] = A/xA$ is a 
polynomial ring in one variable over the field $\kappa(D) = D/\m_D$.}
\end{discussion}

Assume Notation~\ref{not1.3}  and let $\U$ be a finite subset of $Q(D)$.  The ring $\O_\U = \bigcap_{\alpha \in \U}\alpha$
is the irredundant intersection of the regular local rings in $\U$ that are minimal with respect to inclusion. 
Therefore in Corollary~\ref{thm3.1}, we consider finite subsets $\U$ of $Q(D)$ for which the
intersection $\bigcap_{\alpha \in \U}\alpha$ is irredundant.

\begin{corollary}  \label{thm3.1}
Let $\U = \{\alpha_1, \ldots, \alpha_n \}$ be a finite subset of $Q(D)$ such that there are no inclusion 
relations among the $\alpha_i$. Then $\O_\U = \bigcap_{i=1}^n \alpha_i$ is a regular Noetherian domain 
with precisely $n$ maximal ideals. The maximal ideals of $\O_\U$ may be labeled as $\{\m_i\}_{i=1}^n$ 
such that $(\O_\U)_{\m_i} = \alpha_i$ for each $i \in \{1, \ldots, n\}$.   Therefore $\U$ defines an
irredundant essential representation of $\O_\U$.
\end{corollary} 

\begin{pf}  By Theorem~\ref{minimal2},  there exists a complete $\m_D$-primary 
ideal $J$ such that $\Proj D[Jt]$ is a nonsingular projective surface over $ D$ and each of the 
$\alpha_i$ is a closed  point on the surface $\Proj D[Jt]$.  By homogeneous prime avoidance, there exists an
affine component  $A$ of $\Proj D[Jt]$ such that each $\alpha_i$ is a localization of $A$.  
By Remark~\ref{rm2.5} and Theorem~\ref{prop3.7},  $ \U$ defines an irredundant 
 essential representation of $\O_\U$.
If $\m_i$ is the center of $R_i$ on $\O_\U$,  then $(\O_\U)_{\m_i} = R_i$.  Since
the nonunits of $\O_\U$ are the elements in $\bigcup_{i=1}^n \m_i$,   it follows that 
 $\O_\U$ has precisely $n$ maximal ideals.  This proves  Corollary~\ref{thm3.1}. 
\qed \end{pf}

Example~\ref{ex111} illustrates Corollary~\ref{thm3.1}. 
   
 \begin{example} \label{ex111} {\em 
 Assume Notation~\ref{not1.3}.     
 Let $\alpha := D[y/x]_{(x, y/x)}$ and $\beta  = D[x/y]_{(y, x/y)}$ be the infinitely near points to $D$
 in the $x$-direction and $y$-direction. Define $D^* $ to be the local quadratic transform of $\alpha$ in
 the $y$-direction, and $D^{**} $ to be the local quadratic transform of $\beta$ in the $x$-direction. Thus
 $$
 D^* = D[y/x, x^2/y]_{(y/x, x^2/y)} \quad \text{  and }  \quad D^{**} = D[x/y,  y^2/x]_{(x/y, y^2/x)}
 $$
 Let $V$ denote the order valuation ring, $\ord_D$,   of $D$,  and let $V_\alpha$ and $V_\beta$ denote the order
 valuation rings of $\alpha$ and $\beta$, respectively.  Observe that $D^* \subset V$ 
 with $D^*_{(x^2/y)D^*} = V$,   and $D^{**} \subset V$ with $D^{**}_{(y^2/x)D^{**}} = V$.\footnote{In classical terminology, $D$ is {\it proximate} to both $ D^*$ and $D^{**}$.} We also have
 $V_\alpha = D^*_{(y/x)D^*}$, and $V_\beta = D^{**}_{(x/y)D^{**}}$.
 
  Define  $E = D^*  \cap D^{**}$.    Corollary~\ref{thm3.1}   implies 
  that  $E$ has precisely 2 maximal ideals and that $D^*$
  and $D^{**}$ are localizations of $E$.   We give the following direct proof:  
  
  Let $J = (x^4, x^2y,  xy^2, y^4)D$ denote   the product of the simple complete ideals 
  $$
  (x, y)D, \quad   (x^2, y)D,  \quad (x, y^2) \text{  associated to the infinitely near points  } 
  D, ~ \alpha, ~ \beta,
  $$
  respectively. 
  
  The complete ideal $J$ is saturated in the sense of Definition~\ref{def4.1}. By Proposition~\ref{sat = ns}, $X = \Proj D[Jt]$ is 
  a nonsingular projective model over $ D$ and the points $D^*$ and $D^{**}$ are both on
  this model.  Define 
  $$
  A~ = ~  D[\frac{J}{x^2y + xy^2}] ~ = ~ D[\frac{x^3}{y(x+y)},~  \frac{x}{x+y},~ \frac{y}{x+y}, ~
  \frac{y^3}{x(x+y)}].
  $$
  Then  $ A$ is an affine component of the projective model $X$.   We observe that 
  $A$ is a subring of both $D^*$ and $D^{**}$.  Notice that the $x$-adic and $y$-adic valuation 
  rings of $D$ do not contain either $D^*$ or $D^{**}$. Moreover,   the $(x+y)$-adic
  valuation ring  $W = D_{(x+y)D}$ of $D$ does not contain  either $D^*$ or $D^{**}$.  The transform
  of $x+y$ from $D$ to $D[y/x]$ is computed by setting $y_1 = y/x$.  Then $x+y = x + xy_1 = 
  x(1 + y_1)$.  Hence $1 + y_1$ is the transform of $x + y$,  and $W = D[y/x]_{(1 + y/x)D[y/x]}$
  does not contain $D^* \subset D[y/x]_{(y/x)D[y/x]}$. Similarly,  $W$ does not contain $D^{**}$. 
  
  The height one prime ideals  of $D^*$ that contain $y(x+y)$ are the centers of $V$ and $V_\alpha$ 
  on $D^*$, and the height one prime ideals of $D^{**}$ that contain $x(x+y)$ are the centers of $V$ 
  and $V_\beta$.  That $A \subset D^*$ and $A \subset D^{**}$ follows because
$$
x^3V \subset y(x+y)V, \quad y^3V \subset x(x+y)V, \quad xV = (x+y)V, \quad  yV = (x+y)V,
 $$ 
 $$
 x^3V_\alpha  =  y(x+y)V_\alpha, \quad y^3V_\alpha \subset x(x+y)V_\alpha, \quad 
 xV_\alpha = (x+y)V_\alpha, \quad  yV_\alpha \subset  (x+y)V_\alpha,
 $$
 and 
  $$
 x^3V_\beta  \subset  y(x+y)V_\beta, \quad y^3V_\beta =  x(x+y)V_\beta, \quad 
 xV_\beta \subset  (x+y)V_\beta, \quad  yV_\beta  =  (x+y)V_\beta.
 $$
The center  of $D^*$ on $A$  is $(x, y, \frac{x^3}{y(x+y)},  \frac{y^3}{x(x+y)},   \frac{y}{x+y})A$,
while the center of $D^{**}$ on $A$ is 
$(x, y, \frac{x^3}{y(x+y)},  \frac{y^3}{x(x+y)},   \frac{x}{x+y})A$. Since 
these are distinct prime ideals, it follows that $E = D^* \cap D^{**}$ has two distinct maximal ideals
and $D^*$ and $D^{**}$ are the localizations of $E$ at these maximal ideals.}
 \end{example}

\begin{remark}  \label{rm1.5} {\em
It can happen that $R$ and $S$ are 2-dimensional regular local rings with the
same quotient field $F$,  and $R \cap S$ is local and is properly contained in both $R$ and $S$. 
Corollary~\ref{thm3.1}  implies this 
cannot happen   if $R$ and $S$ are both overrings of a 2-dimensional regular
local ring $D$.    

The following example is given  in \cite{H1973}:  Let $x$ and $y$ be variables over a field $k$ and 
let  $R = k[x, x^2y]$  localized at the maximal ideal  $(x, x^2y)R$  and let 
$S = k[xy^2, y]  $ localized at the maximal ideal $(xy^2, y)$.  Then $R$ and $S$  both have quotient field 
$k(x, y)$  and both are  subrings of the formal power series ring  $k[[x, y]]$.   Every element in 
$k[[x, y]]$  has a unique expression as an infinite sum of monomials  in $x$ and $y$ 
 with coefficients from $k$.  Every element in $R$ regarded as a formal power series in $k[[x, y]]$  has the property that the  
$x$-degree of each monomial is greater than the $y$-degree, and for any element in 
$S$ the $y$-degree of each monomial  is greater than the $x$-degree. Hence $R \cap S = k$.  
}
\end{remark}

  \section{Noetherian intersections } \label{sec7}

As in Notation~\ref{not1.3},  let $\R(D)$ denote the overrings of $D$ obtained as intersections 
of rings in $Q(D)$.  As an intersection of normal rings, each ring in $\R(D)$
is  a normal
overring of $D$.
 In Lemma~\ref{lemma 2.1} and Theorem~\ref{normal}, we 
  make some general observations about normal overrings of $D$.

\begin{lemma} \label{lemma 2.1} Let $R$ be a normal overring of a $2$-dimensional Noetherian domain.  Then 
\begin{enumerate}[(1)]
\item $\dim R \leq  2$.   If $R$ has a 2-dimensional Noetherian overring, then $\dim R = 2$. 

\item If ${\ff p}$ is a nonzero nonmaximal prime ideal of $R$, then $R_{\ff p}$ is a DVR and $R/{\ff p}$ is a Noetherian domain.
  
\item $R$ is a Krull domain if and only if $R$ is a Noetherian domain. 

\end{enumerate}
\end{lemma} 

\begin{pf} For item~1, as an overring of the $2$-dimensional Noetherian domain $D$,  it follows 
from the dimension inequality \cite[Theorem 15.5]{Mat} that  every finitely generated $D$-algebra overring
of $D$ has dimension at most 2, and this implies 
$\dim R \leq 2$.   
Assume that $R$ has a Noetherian overring $A$
 with $\dim A = 2$. Let $\m$ be a maximal ideal of $A$ with $\hgt \m = 2$. 
  Since $ A$ is Noetherian, a nonzero element of $A$ is contained in only finitely many height 1 
  primes $\p$ of $A$, and   $\p \cap R \ne 0$,   since $A$ is an overring of $R$.  
  There exist infinitely many height 1 prime ideals $\p$ of $A$
  with $\p \subset \m$.  It follows that $(0) \subsetneq \p \cap R \subsetneq \m \cap R$, 
  for some $\p$ and $\dim R \ge 2$.

  For item~2, see \cite[Proposition 2.3]{ONoeth}, and for item~3, see \cite[Theorem 9]{HKrull}. 
\qed \end{pf}

 \begin{theorem} \label{normal} Let $D$ be a $2$-dimensional regular local ring, and let $R$ be a normal overring of   $D$. 

\begin{enumerate}[(1)]

\item  Assume that  each maximal ideal of $R$ has height $2$.  Then $R$ is Noetherian 
$\iff$      $R$ is a flat overring of a finitely generated $D$-subalgebra of $R$.

\item If $R$  is  local, Noetherian and $\dim R \geq 2$, then $R$ is a spot over $ D$, that is,  $R$ is essentially finitely 
generated over $D$ in the sense that $R$ is the localization of a finitely generated $D$-algebra.

\item A $2$-dimensional  normal local  Noetherian overring $T$ of $R$ is a localization of $R$ $\Longleftrightarrow$ 
each height $1$ prime ideal  of $T$ contracts to a height $1$ prime ideal of $R$.

\end{enumerate}
 \end{theorem}
 
 \begin{pf} 
In  item~1,  the $\impliedby$  direction is clear because every ideal in a flat overring is an extended ideal.
To prove $\implies$  assume $R$ is Noetherian and  let $\p$ be a height~1 prime ideal of $R$ and let $\q = \p \cap D$.  Either $\q = \m_D$
or $\hgt \q = 1$.  If $\hgt \q = 1$,  then $D_\q$ is a DVR,  and $D_\q = R_\p$.   A nonzero element in 
the  Noetherian domain  $R$ is contained in only finitely many height 1 primes of $R$.  Hence
there exists only  a finite set, say $\{\p_1, \ldots, \p_n\}$,   of height 1 prime ideals 
of $R$ that contain
 $\m_D$. 
 
 Since each maximal ideal of $R$ has height 2,  the DVRs $V_i = D_{\p_i},  i \in \{1, \ldots, n\}$ are prime
 divisors of the second kind over $D$. 
 Hence there exist elements $a_i \in R$ such that the image of $a_i$
 in the residue field of $V_i$ is algebraically independent over $D/\m_D$.
 
 Let $A$ denote the integral closure of $D[a_1, \ldots, a_n]$.   A classical result of Rees \cite{Rees}
 implies that   $A$ is a finitely generated $D$-algebra.   Thus $A$ is a 
 normal Noetherian subring of $R$ such that for each height 1 prime $\p$ of $R$, 
 then ht$(\p \cap A) = 1$. 
 
 Let $\m$ be a maximal ideal of $R$ and let $\n = \m \cap A$. If $\hgt \m = 1$, then $\hgt \n = 1$ and 
 $A_\n = R_\m$ is a DVR. If $\hgt \m > 1$, then $A_\n \subseteq R_\m$ are 2-dimensional normal 
 Noetherian local  domains with $R_\m$ dominating $A_\n$. 
 Let  $\q \in \Spec A_\n$ with $\hgt \q = 1$. 
 Since $A_\n$ dominates $D$, $A_\n$ has a rational singularity. 
 By \cite[Proposition 17.1]{L1969}, $A_\n$ has a 
finite divisor class group. Hence $\q$ is the radical of a principal ideal of $A_\n$. 
It follows that $\q R_\m$ is contained in a height 1 prime of $R_\m$. 
Therefore the set of essential valuation rings for $A_\n$ is the same as the set of 
essential valuation rings for $R_\m$,  and thus $A_\n = R_\m$.  
By Theorem~2 of \cite{Ric}, $R$ is a flat overring of $A$.  
Item~2 follows from item~1. 

For item~3, it is clear that if $T$ is a localization of $R$, then each height $1$ prime ideal of $T$ contracts to a height $1$ prime ideal of $R$.  Suppose $T$ is not a localization of $R$.  
By replacing $R$ with $R_{{\ff m}_T\cap R}$ we may assume without loss of generality that $R$ is a 
 normal local  ring with $R \subseteq T$.  
 By item 1, there exist  $x_1,\ldots,x_n \in T$ such that $T$ is a localization of $D[x_1,\ldots,x_n]$. 
 Hence $T$ is a localization of $A := R[x_1,\ldots,x_n]$.  
 Since $R$ is integrally closed,   $R$ is  integrally closed in  $A$.  Peskine's version of Zariski's Main Theorem \cite[Proposition 13.4, p.~174]{Pes} implies that  there is a height one prime ideal ${\ff p}$ of $T$ such that ${\ff p} \cap R = {\ff m}_{T} \cap R$.  Since $\dim T = 2$, Lemma~\ref{lemma 2.1}.2 implies ${\ff m}_T \cap R$ is a height $2$ prime ideal of $R$. Therefore if $T$ is not a localization of $R$, there is a height $1$ prime ideal of $T$ that does not contract to a height $1$ prime ideal of $R$.  
 This 
 proves item~3.  
\qed \end{pf}

Let $R$ be a normal Noetherian overring of $D$. By a desingularization of $\Spec R$ we mean a desingularization $Y$ of the model $X= \{R_{\ff p}:{\ff p} \in \Spec R\}$  over $R$; i.e., $Y$ is a nonsingular projective model over $R$ that dominates $X$.

\begin{theorem} \label{desing} Let $R$ be a normal Noetherian overring of $D$ for which every maximal ideal has height $2$. Then 
\begin{enumerate}[(1)]
\item There exists a desingularization $Y$ of $\Spec R$ such that $Y$ is a subset of  
 a nonsingular projective model 
$X = \Proj D[Lt]$  over $D$.
\item $\Spec R$ has finitely many singularities.
\item Each desingularization of $\Spec R$ is a product of quadratic transformations.
\end{enumerate}
\end{theorem}

\begin{pf} For item 1, by Theorem~\ref{normal}.1 there exists a 
finitely generated $D$-subalgebra $A$ of $R$ such that $R$ is flat over $A$.   We may assume that
$a_1, \ldots , a_n,  b$ are nonzero elements in $D$ such 
that $A = D[a_1/b,  \ldots , a_n/b]$.  For each maximal ideal $\m$ of $R$, the local ring $R_\m$ is a flat 
overring of $A$ and hence $A_{\m \cap A} = R_\m$.  

Let $I = (a_1, \ldots , a_n, b)D$.  Then each $R_\m$ is on the model $\Proj D[It]$. Let $J$ denote the
integral closure  of $I$. Since $R$ is integrally closed,  each $R_\m$ is on the 
normal model $\Proj D[Jt]$.  
Let $L$ be the saturation of $J$.  Then $X =  \Proj D[Lt]$ is a nonsingular projective model  over $\Spec D$
 that dominates $\Proj D[Jt]$, and the map $f: X \to \Proj D[Jt]$ is a 
desingularization  of  $\Proj D[Jt]$.  

Since $\Spec R \subset \Proj D[Jt]$, the inverse image in $X$   of  $\Spec R$ with respect to the map $f$
is a  desingularization of $\Spec R$. This proves item~1.


To prove item 2,  by
Theorem~\ref{normal}.2,  each localization of $R$ at a maximal ideal is a normal spot over $D$. 
By \cite[Proposition, p.~160]{Lip2}, each normal spot  over $D$ is analytically normal. 
That $\Spec R$ has finitely many singularities follows now from    \cite[Theorem, p.~151]{Lip2} .


To prove item 3, we may assume 
that $\Spec R$ has singularities. 
 Since  each localization of $R$ at a maximal ideal is a spot over $D$ by Theorem~\ref{normal}.2, it follows  from   
   \cite[Proposition~1.2]{L1969} that  each normal Noetherian local overring  of $D$  has a rational singularity. Thus each of the finitely many singularities of $\Spec R$  is a rational singularity. 
A result of  Lipman  \cite[Theorem~4.1]{L1969} implies that any desingularization of $\Spec R$ is a product of quadratic transformations.
\qed \end{pf}

With Theorem~\ref{desing}, we obtain a characterization of the Noetherian rings in $\R(D)$.  

\begin{theorem} \label{cor3.6} 
Assume Notation~\ref{not1.3}. 
The following are equivalent for an overring $R$ of $D$.
\begin{enumerate}[(1)]
\item $R$ is a Noetherian domain   in $\R(D)$.
\item $R$ is a  normal Noetherian domain  for which   every maximal ideal  has height~$2$.
\item  There exists a nonsingular projective model $X$ over $D$ and a subset $\U$ of
the closed points of $X$ such that $R = \O_{\U}$.

 
 
\end{enumerate}
\end{theorem} 

\begin{pf}
(1) $\Longrightarrow$ (2): 
Since $R \in \R(D)$, there exists a subset $\U$ of $Q(D)$ such that 
$R = \O_\U = \bigcap_{\alpha \in \U} \alpha$.  Let $\m_{\alpha}$  denote the maximal
ideal of $\alpha$. Lemma~\ref{lemma 2.1}.1  implies that ht$(\m_{\alpha} \cap R) = 2$
for each $\alpha \in \U$.  
Therefore 
$$
R  ~  \subseteq \bigcap_{\alpha \in \U}R_{\m_{\alpha} \cap R} ~\subseteq \bigcap_{\alpha \in \U} \alpha ~= R.
$$
Suppose there exists a height 1 maximal ideal $\p$ of $R$.
Then $\p$ is invertible and 
$R \subsetneq \p^{-1}$.   But $\p^{-1} \subseteq R_\m$ for each maximal ideal $\m$ of $R$
with $\hgt \m = 2$,  and hence  $\p^{-1} \subseteq R$,  a contradiction.  
Therefore if $R \in \R(D)$ is Noetherian, then every maximal ideal of $R$ has 
height 2.\footnote{It is shown in  \cite{HLO}  that there exist non-Noetherian $R \in \R(D)$ that have 
maximal ideals of height 1.}

(2) $\Longrightarrow$ (3): By Theorem~\ref{normal}.1,  there exists a 
 projective model $\Proj D[Jt]$  such that $\{R_{\ff p}:{\ff p} \in \Spec R\}   \subset \Proj D[Jt]$.  By Theorem~\ref{thm2.6}, 
 there is a desingularization $X$ of $\Proj D[Jt]$.
 Let $\U$ be the closed points in $X$ that contain $R$. 
 Then $R = \O_{\U}$.



(3) $\Longrightarrow$ (1): Each of the local  rings in $\U$ is an intersection of exceptional prime
 divisors of $X$  and prime divisors of the first kind. Since there are only finitely many exceptional prime divisors of $X$ and the set of prime divisors of $D$ of the first kind has finite character, it follows that $R =\O_S$ is a finite character intersection  of DVRs, hence a Krull domain. Lemma~\ref{lemma 2.1} implies that as a Krull overring of a $2$-dimensional Noetherian domain, $R$ is a normal Noetherian domain.  
\qed \end{pf}

\begin{corollary}  \label{cor3.4}
Let $X$ be a nonsingular projective model over $D$,  and let $L$ be a saturated complete 
ideal of $D$ such that $X = \Proj D[Lt]$.  
\begin{enumerate}[(1)]
\item
There exist only finitely many 
 local  domains $R$  that are not regular and have 
the form $R = \O_{\U}$,  where $\U$ is a subset of the
closed points of $X$. 
\item
Each  $R$   is normal Noetherian with $\dim R = 2$.
\item
 Each  $R$  in item~1 is on a normal projective model $N = \Proj D[Jt]$ 
 over $D$,  where $J$ is a complete ideal of $D$ that divides $L$.
 \end{enumerate}
\end{corollary} 

\begin{pf}
If $R = \O_{\U}$,  where $\U$ is a subset of the
closed points of $X$,   then Theorem~\ref{cor3.4} implies that 
$R$   is normal Noetherian with $\dim R = 2$,  and 
$R$ is a point on a normal projective model that is dominated by $X$. 
Every  normal projective model dominated by $X$  has the form $\Proj D[Jt]$, where $J$ is a complete 
ideal that is the product of a subset of the simple complete factors of $L$.  Two subsets with the same
simple  complete factors define the same model. Hence there exist only finitely many normal projective
models over $D$ that are dominated by $X$.  Theorem~\ref{desing} implies that each of these normal 
projective models has only finitely many singular points.
Therefore there are only finitely many $R$ of this form that are not regular. 
\qed \end{pf}

  
\begin{example}  \label{ex1.6} {\em 
 Assume Notation~\ref{not1.3}.   Let $R = D[y^2/x]_{(x, y, y^2/x)D[y^2/x]}$. 
  Then 
 the maximal ideal $\m_R$ of $R$ is  $(x, y,  y^2/x)R$.  Since $\m_R$ requires 3 generators, $R$ is not regular.  $R$ is a normal domain,  and 
 is  called an {\it  ordinary double point singularity}.   We give an explicit representation of $R$ as an intersection of rings in $Q(D)$.
 
  We show that  the quadratic transform $\Proj R[\m_Rt]$ is a nonsingular model over $ R$.  As a normal domain, $R$ is an intersection of its minimal valuation overrings. Moreover, each minimal valuation overring of $R$ dominates a closed point in $\Proj R[\m_Rt]$. We use this fact in describing the closed points in $\Proj R[\m_R t]$.

Let $V$ be a minimal valuation overring of $R$. Then $V$ dominates $R$, and $\m_RV$ is either
$xV$, or $yV$ or $(y^2/x)V$.

If $\m_RV = xV$, then $y/x \in V$ and $R[y/x] = D[y/x] \subset V$.  The affine component $D[y/x]$ of 
$\Proj R[\m_Rt]$ is nonsingular,  and   $V$ is centered on a 
height 2 maximal ideal of $D[y/x]$.  
Since $y^2/x \in yD[y/x]$ and $y \in xD[y/x]$, it follows that $xD[y/x] \cap R = (x, y, y^2/x)R$.

Every maximal ideal  $\p$ of  $D[y/x]$ of height 2 contains $x$,
 and   the map 
$\Spec D[y/x] \to \Spec R$   maps each maximal ideal $\p$ of 
height 2 of $D[y/x]$  to  
$\p \cap  R  =  \m_R$. Therefore    all the elements of
$Q_1(D)$ other than $\alpha  = D[x/y]_{(y, x/y)D[x/y]}$ dominate $R$, and are dominated
by a minimal valuation overring of $R$. 

If $\m_RV = yV$,  then both $x/y$ and $y/x$ are in $V$.  Hence the affine component of $\Proj R[\m_Rt]$
 obtained by dividing by $y$ gives nothing  that is not already obtained in the affine component 
 dividing by $x$.

If $\m_RV = (y^2/x)V$,  the affine component obtained by dividing by $y^2/x$ is $R[x/y] = D[x/y, y^2/x]$ and is 
nonsingular.  If $xV \ne yV$, then $V$ is centered on the maximal ideal  $(x/y, y^2/x)D[x/y, y^2/x]$.  The localization
at this maximal ideal is the point 
 $\beta   =  \alpha[y^2/x]_{(x/y, y^2/x)\alpha[y^2/x]}  \in Q_2(D)$.   
 Notice that $\beta$ dominates $\alpha \in Q_1(D)$.

Let $A = \bigcap\{D[y/x]_\p ~|~ \hgt \p = 2 \}$.
Then  $A$ is the intersection of the elements in $Q_1(D)$ other than $\alpha$.  
The centers on $\Proj R[\m_Rt]$  of the minimal valuation overrings of $R$ are  the localizations of 
$A$ at its maximal ideals and $\beta$.  Therefore $\Proj R[\m_Rt]$ is nonsingular,  and 
 $R = A \cap \beta$. 

Since $R$ is proper subring of $A$, the point $\beta$ is  irredundant in this representation. 
Let $\gamma$ be a point in $Q_1(D) \setminus \{\alpha\}$.  Theorem~\ref{prop4.2} implies that
there exists an element $f \not\in D$ such that   $f$ is in each point
of  $Q_1(D) \setminus \{\gamma\}$.  Since $\beta$ dominates $\alpha$,  $f \in \beta$. 
Since $R \cap \alpha = D$,  it follows that $f \not\in R$.  Therefore the representation 
$R = A \cap \beta$ is irredundant.

$R$ is the unique singular point of the normal projective model $\Proj D[xt, y^2t]$.  The 
quadratic transform $\Proj R[\m_Rt]$ described above is in the nonsingular model 
$\Proj D[Jt]$,  where $J = (x^2, xy, y^3)D$. 

The saturated complete ideal $J = (x, y)(x, y^2)R$ has two complete  simple factors.
Corollary~\ref{cor3.4}   implies that $R$ is the unique non-regular local domain 
of the form
 $ \O_{\U}$,  where $\U$ is a subset of the
closed points of $\Proj D[Jt]$.  }
\end{example}

   



\section{Irredundant representations}   \label{sec8}


Each ring $R$ in $\R(D)$ is an intersection of rings in $Q(D)$. In this section we consider settings in which $R$ can be represented by an irredundant intersection
of rings from $Q(D)$; i.e., there is a set $\U$ in $Q(D)$ such that $R = \O_\U$ but $R \subsetneq \O_{\U \setminus \{\alpha\}}$ for each $\alpha \in \U$.

\begin{definition} {\em  A subset $\U$ of $Q(D)$ is said to be {\it complete} 
if $\O_\U = \bigcap_{R \in \U}R = D$.  For a point $\alpha \in Q(D)$, a subset $\U$ of $Q(D)$ of points
that dominate $\alpha$ is said to be {\it complete over } $\alpha$ if $\O_\U = \alpha$. } 
\end{definition}

\begin{remark} \label{rem5.3} {\em 
Let  $X = \Proj D[Jt]$ be a normal  projective model over $ D$. Then 
the set $\U$ of closed points of $X$ is complete since every minimal valuation ring $V $ 
is  centered on a closed point of $X$ and $D$ is the intersection of the minimal  valuation overrings of $D$.  
It is natural to ask if a proper subset of $\U$ 
can  be complete. This is equivalent to asking if the representation $D = \bigcap_{\alpha \in \U}\alpha = D$
is irredundant.  }
\end{remark}

\begin{theorem} \label{prop4.2}
 Assume notation as in Discussion~\ref{dis3.3}. 
 \begin{enumerate}[(1)]
  \item   If  $\U = Q_1(D)$,  then $\U$  is complete, and the 
  representation $D = \bigcap_{R \in Q_1(D)}R$ is
 irredundant.
 \item  Assume $\U$ is a nonempty  proper subset of $Q_1(D)$,  then
 \begin{enumerate}[(1)]
 \item
 $\O_\U$ is a flat extension  of a regular finitely generated $D$-subalgebra of $F$.
 \item 
 $B := \O_\U$ is a regular Noetherian domain and the representation $B = \bigcap_{\alpha \in \U} \alpha$
 is an irredundant essential representation of $B$.  
 \end{enumerate}  
 \end{enumerate}
\end{theorem}

 \begin{pf}
  The set $\U = Q_1(D)$ is complete by Remark~\ref{rem5.3}.
   To see that  the representation $D = \bigcap_{R \in Q_1(D)}R$ is
 irredundant,   we use as in Remark~\ref{rem3.3} that  the points of $Q_1(D)$  are in a natural one-to-one correspondence with 
 the maximal homogeneous relevant\footnote{A homogeneous  prime ideal  of the graded domain
 $D[xt, yt]$  is said to be {\it relevant} if it does not contain the maximal graded ideal 
 $(x, y, xt, yt)D[xt, yt]$.} prime ideals of $D[xt, yt]$.
  The relevant homogeneous maximal ideals of $D[xt, yt]$ all contain $\m_D$, and
 $$
 D[xt, yt]/\m_DD[xt, yt]   ~ \cong   ~ \kappa(D)[\overline x, \overline y],  ~\text{  where }~ 
 xt \mapsto \overline x ~ \text{  and  } ~ yt \mapsto   \overline y,
 $$
  and $\overline x, \overline y$ are algebraically independent over $\kappa(D)$. Hence to each 
  point $\gamma \in Q_1(D)$ there corresponds an irreducible homogeneous polynomial 
  $f(xt, yt) \in D[xt, yt]$ such that the image $\overline f$ of $f$ in $\kappa(D)[\overline x, \overline y]$
  is irreducible.  Then $\deg f = \deg \overline f = d$ for some positive integer $d$,
  and the degree zero component $A  = D[\frac{(x, y)^d}{f(x, y)}]$   
   of the $\mathbb Z$-graded ring $ D[xt,  yt][1/f]$ is an affine component of $\Proj D[xt, yt]$.
   Then $A$ is a  regular finitely generated $D$-subalgebra of $F$.

   Let $\U = Q_1(D) \setminus \{\gamma\}$.
   Each of the points of $\U$ is a localization of  $A$.    Therefore 
   $A \subseteq \O_\U$.  Since $A$   properly
  contains $D$,  it follows that $\gamma$ is irredundant in the representation 
  $D = \bigcap_{R \in Q_1(D)}R$.\footnote{The local ring $\gamma$ does not contain $A$.}
  Since this is true for each $\gamma \in Q_1(D)$, the representation 
   $D = \bigcap_{R \in Q_1(D)}R$ is irredundant.  This proves item 1.
   
   For item 2,  assume that $\U$ is a nonempty proper subset of $Q_1(D)$.   Then $\U$ is a subset of 
   $Q_1(D) \setminus \{\gamma\}$,  for some $\gamma \in Q_1(D)$. 
   Hence $B = \O_\U$ is a flat extension of the regular finitely generated $D$-algebra    
    $A  = D[\frac{(x, y)^d}{f(x, y)}]$  and the representation $B = \bigcap_{\alpha \in \U} \alpha$
 is an irredundant essential representation of $B$.   
 \qed \end{pf}

 Let $\alpha \in Q(D)$. 
 Applying Theorem~\ref{prop4.2} to the 2-dimensional regular local domain $\alpha$ 
 gives Corollary~\ref{Q_1}.


\begin{corollary} \label{Q_1}  Let  $\alpha \in Q(D)$.  Let $\U$ be a  subset of $Q_1(\alpha)$. If $\U = Q_1(\alpha)$, then $\O_\U = \alpha$. Otherwise, if $\U$ is a proper subset of $Q_1(\alpha)$,  then
\begin{enumerate}[(1)]
\item  $\O_\U$ is a  regular Noetherian domain such that $\U$ is the set of localizations at a maximal ideal is in $\U$.

\item $\O_\U$ is a flat extension of a regular finitely generated $D$-subalgebra of $F$. 

\item The representation of $B = \O_\U =  \bigcap_{\alpha \in \U}\alpha$ is an irredundant essential 
representation. 
\end{enumerate} 
\end{corollary}

The rest of the section is devoted to the case where $D$ is Henselian. We first establish a lemma that applies to the Henselian case but whose hypotheses can hold in more general settings for specific choices of height one prime ideals of $D$.


  
\begin{lemma}  \label{prop4.5}
Assume Notation~\ref{not1.3}.
Let $\gamma \in Q_1(D)$ and let $\p$ be a height one prime of $D$ such that $\gamma \subset D_\p$. 
\begin{enumerate}[(1)]
\item  
Then  
$\gamma/(\p D_\p \cap \gamma)$ is a local quadratic transform of $D/\p$.  
\end{enumerate}
Suppose in addition that  the integral closure of $D/{\ff p}$ is local. 
\begin{enumerate}[(1)]
\item[{\em (2)}] 
The ring  $\gamma$ is the only point of $Q_1(D)$ such that $\gamma \subset D_\p$.
\item[{\em (3)}] 
For each integer $n \ge 2$,  there exists a unique 
point $\gamma_n \in Q_n(D)$ such that $\gamma_n \subset D_\p$.  Then 
$V = \bigcup _{n=1}^\infty\gamma_n $ is
the rank 2 valuation overring of $D$ that is the composite of $D_\p$ with the integral closure of $D/\p$.

\end{enumerate}
    \end{lemma}  
    
\begin{pf} 
For item~1,  since  $\gamma \subset D_\p$,   the  canonical surjective map $D_\p \to D_\p/\p D_\p$ restricts to a
surjective map $\gamma \to \gamma/(\p D_\p \cap \gamma)$.  
Since $\gamma$ is a local quadratic transform of $D$,  the universal property of blowing up
\cite[Prop.~7.14 and Cor.~7.15, pp.~164--165]{H} implies  that 
the induced map $ D/\p   \to \gamma/(\p D_\p \cap \gamma)$ is a local quadratic transform.

For item~2,  let  $R = D/\p$.   Then $R$ is a Noetherian local domain with $\dim R = 1$.  Since the 
integral closure of $R$ is local, the integral closure of $R$ is the unique valuation overring of $R$ 
dominating $R$.  It follows that every overring of $R$ is local. 

 If $\gamma' \in Q_1(D)$ with $\gamma' \ne \gamma$ and $\gamma' \subset D_\p$,
then  Corollary~\ref{thm3.1}   implies that 
$A = \gamma \cap \gamma'$ has two maximal ideals   and both of the 
maximal ideals of $A$ contain $\p D_\p  \cap A$. This implies that $A/(\p D_\p \cap A)$ is an
overring of $R$ that is not local,  a contradiction. Hence $\gamma$ is the unique point in $Q_1(D)$
with $\gamma \subset D_\p$.  

A similar argument proves for each $n \ge 2$  that there exists a
unique point $\gamma_n \in Q_n(D)$ such that $\gamma_n \subset D_\p$.  This 
proves item~3.
\qed \end{pf}  

If $D$ is Henselian, then the integral closure of $D/\p$ is local for each height 1 prime $\p$ of $D$. Thus the statements in items~2 and 3 apply to
$D_\p$.  We use this observation in the proof of Theorem~\ref{thm5.10}.

\begin{theorem} \label{thm5.10}
Assume that $D$ is Henselian and  $\U$ is a set of pairwise incomparable rings in $Q(D)$. 
Then the representation $\O_\U = \bigcap_{\alpha \in \U} \alpha$
 is irredundant. 
\end{theorem}

\begin{pf}  We first make a couple of reductions to simplify  the proof. 
Let $\U^*$ be the union of $\U$ with the  set of rings $\alpha \in Q(D)$ such that $\alpha$ is minimal with respect to not containing any of the rings in $\U$. Then $D = \O_{\U^*}$ and to prove that $\U$ is an irredundant representation of $\O_\U$ it suffices to prove that $\U^*$ is an irredundant representation of $D$.  Thus we may assume that $\U  = \U^*$   and hence that $\U$ is complete.

 Let $\gamma \in \U$. It suffices to show that $D \subsetneq \O_{\U \setminus \{\gamma\}}$, and to prove this it suffices to show that $D \subsetneq \O_{\U'}$, where  $\U'$ is the set of rings in $Q(D)$ minimal with respect to not containing $\gamma$.  By Theorem~\ref{minimal2}.1 we may assume without loss of generality that $\U$ is the set of closed points in a projective nonsingular model over $D$.  

  Let $J$ be a complete ideal of $D$ such that $X = \Proj D[Jt]$.  
Let $\E(D)$ denote the set of essential valuation rings of $D$. For each $\alpha \in \U$, the set $\E(\alpha)$ of essential
valuation rings for $\alpha$ is the union of a subset of $\E(D)$ with a subset of $\Rees J$. 
Since $\Rees J$ is a finite set and $\E(\alpha)$ is an infinite set,
for each $\alpha \in \U$,  there exists a height~1 prime
$\p_{\alpha}$ of $D$ such that  $D_{\p_{\alpha}} \in \E(\alpha)$.

Since $D$ is Henselian,  the integral closure of $D/\p_{\alpha}$ is local.  
Since there are no inclusion relations among the points in $\U$,   Lemma~\ref{prop4.5}.3 
implies that $\alpha$ is the unique point of $\U$ contained in $D_{\p_{\alpha}}$.  

The set $\S : = \E(D) \setminus (\{D_{\p_{\alpha}}:\alpha \in \U\} \cup \Rees J)$ includes the set $\E(\beta)$ for each 
$\beta \in \U$ with $\beta \ne \alpha$.
Therefore $\bigcap\{V ~|~ V \in \S \} ~ \subseteq ~  \bigcap\{\beta ~|~ \beta \in \U, ~ \beta \ne \alpha \}$.

Since the set $\S$ has finite character in the sense that a nonzero element of $F$ is a unit in all but finitely many of the elements in $\S$ and since $\{D_{\p_{\alpha}}\} \not\subseteq \S$,  it follows that 
$D ~ \subsetneq     ~ \bigcap\{V ~|~ V \in \S \}$.   Since $\alpha$ is an arbitrary element in $\U$,
the representation $D = \O_\U$ is irredundant.
\qed \end{pf} 

\begin{corollary} \label{Henselian cor} Assume that $D$ is Henselian. Let $R$ be a Noetherian normal overring of $D$ such that every maximal ideal has height $2$.  Then $R$ is an irredundant intersection of  
 the rings in $Q(D)$ that are minimal with respect to containing $R$.  
\end{corollary}

\begin{pf}
Apply Theorem~\ref{cor3.6} and  Theorem~\ref{thm5.10}. 
\qed \end{pf}







\medskip

{\it Acknowledgment.}  This paper is an outgrowth of a project  the authors are working on with Alan Loper  in \cite{HLO}. We thank Alan for  helpful conversations and for motivating us to consider the topics in this paper. We also thank Joe Lipman for several helpful comments on an earlier draft of this paper.

\medskip


\end{document}